\documentstyle[12pt]{article}   

\title{Dipper-Donkin algebra as global symmetry of quantum chains.} 
\author{     
Suemi Rodr\'{\i}guez-Romo\\Centre of Theoretical Research\\
Universidad Nacional Aut\'onoma de M\'exico, Cuautitl\'an\\
Apdo. Postal 142, Cuautitl\'an Izcalli\\
Estado de M\'exico, 54740 M\'exico.$^{\ast}$}         
\date{}  
\begin{document}   
\maketitle
\renewcommand{\thefootnote}{\fnsymbol{footnote}}
\setcounter{footnote}{-1}
\footnote{$\hspace*{-6mm}^{\ast}$
e-mail: suemi@servidor.unam.mx }
\renewcommand{\thefootnote}{\arabic{footnote}}
\baselineskip0.6cm
{\small {\bf Abstract.} 
We analyze the role of $GL_2$, a quantum group constructed by 
Dipper and Donkin \cite{dido}, as a global symmetry for quantum 
chains, and show the way to construct all possible Hamiltonians 
for four states quantum chains with $GL_2$ global symmetry. In 
doing this, we search all inner actions of $GL_2$ on the Clifford 
algebra ${\it C}(1,3)$ and show them. We also introduce
the corresponding operator algebras, invariants and Hamiltonians, 
explicitly.}
\newpage
\section{Introduction.}
In the last few years, quasi triangular Hopf algebras or quantum groups 
have attracted a lot of attention from physicists. One of the most 
interesting features is that such structures can be related to underlying 
symmetries on spaces where the coordinates are noncommutative. Promoting 
these coordinates to functions, it has been shown that it is possible to 
write down an action for such fields that, when added to the action of a 
commuting field, has a symmetry resembling supersymmetry. Quantum groups 
can help us to understand the transformations on such fields 
and the action invariances.\

Symmetry has always played an important role in theoretical physics 
in helping to reduce a problem with many variables to a more tractable
size. The basis of the method built on Bethe Ansatz is to diagonalize the 
Hamiltonian along with an infinite set of constants of motion. In some
cases the ocurrence of this infinite set of constants of motion is related 
to the appearance of a new kind of symmetry, the quantum group symmetry.
This nourishes the hope that, by relaxing the demands usually made
on the structure of a symmetry group, and allowing the wider class of 
quantum groups, one can benefit from symmetry considerations in new
situations, where a symmetry in the traditional sense is simply not 
present.\

We model physical systems where the variables at the lattice sites
take values so that the operators acting on them are matrices
of dimension four by four and complex entries.\

For some subclass of conformal integrable systems; well-know examples 
are given by the minimal models, the WZW models and Liouville-Toda 
theory, the underlying symmetries are indeed known to be given
by quantum groups. However, in spite of extensive studies, our 
understanding of the quantum group symmetry in these theories is 
still somewhat incomplete. We also think that there are realizations 
of the quantum group symmetry in nature.\

The adjoint inner actions studied here are also called 
spectrum generating quantum group. So, we are classifying all 
possible roles of $GL_2$ as a spectrum generating algebra
for ${\it C}(1,3)$. These inner actions have been used as a gauge 
transformation $H\rightarrow B$ in quantum group gauge theory. Here 
$B$ is the space-time algebra and $H$ the coordinate ring of the 
gauge group.\

Quantum $GL_n$ are unique because they are related with $q$-Schur 
algebras \cite{Ja}, hence Hecke algebras and the representation theory
of finite general linear groups. In fact the representations of quantum 
$GL_n$ provide the connection between the classical theory of
polynomial representations of infinite general linear groups and the 
representation theory of finite general linear groups in the 
non-describing characteristic case; obviously if we take $q$ to be 1, we 
are in the classical case. In this limit the representation theory of 
Dipper-Donkin quantum groups is equivalent to the representation theory of 
$q$-Schur algebras. From this follows that the importance of this paper is 
to gain an understanding of the $q$-Schur algebra as a global symmetry for
quantum chains. We remark that the Dipper-Donkin quantum group are
not special cases of the well known Manin's construction. There is one
fundamental difference: {\it the Dipper-Donkin quantum determinant is in 
general not central}.\

Weyl and Clifford algebras are at the heart of quantum physics. The
most useful of them are those endowed with definite transformation properties
under the action of some symmetry group. The idea that quantum groups
could generalize Lie groups in describing symmetries of quantum physical
systems has attracted much interest in the past decade.\

In this paper we study the inner action of the Dipper-Donkin
quantum group on the ${\it C}(1,3)$ algebra, namely the algebra 
generated by the Dirac matrices, as a testing ground for applications
of quantum group symmetries. We search for the corresponding operator 
and invariant algebras in order to have additional information to propose 
Hamiltonians for quantum chains with this global symmetry. Interpreting 
the quantum group as a gauge group, one would consider only the invariant 
elements as observables. The rest of the algebra would then be an algebra
of unobservable fields, whose function in the theory is to describe
operations changing the superselection sector (creating charge). 
Besides, we are interested in some fundamental questions. Can a quantum 
chain have global symmetry given by a quantum group with no 
central but group-like determinant? and what would it be the meaning of 
this?\

We address here the first question and study (as a particular 
case) four states quantum chains. We are able of showing all possible, 
non-trivial Hamiltonians for this system, with Dipper-Donkin global 
symmetry.\

Having discovered all the Hamiltonians which are invariant under 
the Dipper-Donkin quantum group for four states quantum chains, we 
think of these, as systems whose energy eigenstates organize into 
$GL_2$ multiplets, with no energy splitting among members of the 
same multiplet. This is done in spite of the fact that the $GL_2$ 
determinat is group-like as it should, but not central to the algebra, 
as it is the case.\

\section{Dipper-Donkin algebra.}
\baselineskip0.6cm

The algebraic structure of Dipper-Donkin quantization $GL_2$ 
\cite{dido} is generated by four elements $c_{ij}$, 
$1\leq i,j\leq 2$ with relations which can be presented by the 
following diagram.\

\vspace*{1.5cm}
\hspace*{2cm}
\begin{picture}(0,80)
\put(115.29,2.39){\makebox(0,0)[cc]{$c_{21}$}}
\put(141,70.5){\makebox(0,0)[cc]{$d$}}
\put(168.32,2.39){\makebox(0,0)[cc]{$c_{22}$}}
\put(168.32,50.46){\makebox(0,0)[cc]{$c_{12}$}}
\put(115.29,50.46){\makebox(0,0)[cc]{$c_{11}$}}
\put(117.50,7){\vector(1,1){40}}
\put(115.5,7){\vector(0,1){40}}
\put(168.5,7){\vector(0,1){40}}
\put(115.7,7){\vector(1,2){28}}
\put(147,68){\vector(3,-2){18}}
\put(120.5,3){\line(1,0){40}}
\put(120.5,51.0){\line(1,0){40}}
\put(115.5,55){\line(5,3){18}}
\put(146,63){\line(2,-5){23}}
\put(123.5,42.0){\makebox(0,0)[cc]{$\cdot $}}
\put(131.5,34.0){\makebox(0,0)[cc]{$\cdot $}}
\put(139.50,26.0){\makebox(0,0)[cc]{$\cdot $}}
\put(147.5,18.0){\makebox(0,0)[cc]{$\cdot $}}
\put(156.5,10.0){\makebox(0,0)[cc]{$\cdot $}}
\end{picture}
\vspace*{.5cm}
\label{fig1}
\newline
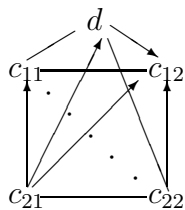
\centerline{Figure 1. Diagramatic representation of}
\centerline{ Dipper-Donkin algebra. 1.
}\\

\vspace{1cm}
Here we denote by arrows $x\rightarrow y$ the ``quantum spinors'' ( or 
generators of the quantum plane \cite{Man}) $xy=qyx$. By 
straigth line $x-y$ we denote the ``classical spinors'' $xy=yx$ and by 
dots $x\cdot\cdot\cdot y$ a classical spinor with a nontrivial 
{\it perturbation} \cite{suemi}, $xy-yx=p$ being
$p$=$(q-1)c_{12}c_{21}$.\

In this algebra the quantum determinant $d=c_{11}c_{22}-c_{12}c_{21}$ 
is noncentral and group-like. This, in contrast with Manin's 
approach \cite{Man}. A group-like element $d$, in a Hopf
algebra, is such that $\Delta d=d\otimes d$ and $\epsilon (d)=1$. 
In any Hopf algebra 
every group-like element is invertible, therefore the quantum $GL_2$ 
includes the formal inverse $d^{-1}$.\

The coalgebra structure is defined in the standard way for all
quantizations and the antipode $S$ is given in reference [1].\

As we know, the Clifford algebra ${\it C}(1,3)$ is generated by the vectors 
$\gamma_{\mu}$, $\mu=0,1,2,3$ with relations defined by the form 
$g_{\mu\nu}$=$diag(1, -1, -1, -1)$, as follows:
$$
\gamma_{\mu}\gamma_{\nu}=g_{\mu\nu}+\gamma_{\mu\nu},\;\;
\gamma_{\mu\nu}=-\gamma_{\nu\mu},
$$
$$
\gamma_{\rho}\gamma_{\mu\nu}=
g_{\rho\mu}\gamma_{\nu}-g_{\rho\nu}\gamma_{\mu}
+\gamma_{\rho\mu\nu},
$$
$$
\gamma_{\lambda}\gamma_{\mu\nu\rho}=g_{\lambda\mu}\gamma_{\nu\rho}-
g_{\lambda\nu}\gamma_{\mu\rho}+g_{\lambda\rho}\gamma_{\mu\nu}+
\gamma_{\lambda\mu\nu\rho}.
$$
This algebra is isomorphic to the algebra of the $4\times 4$ complex 
matrix and it includes the basis of matrix units reported in reference
\cite {aqg}; $e_{ij}$, $1\leq i,j\leq 4$, among others.\

An action of $GL_2$ on ${\it C}(1,3)$ is uniquely defined by actions of
$c_{ij}$ on the generators of ${\it C}(1,3)$\cite{Co}\cite{Sh};
\begin{equation}
c_{ij}\cdot \gamma_k=
f_{ijk}(\gamma_0, \gamma_1, \gamma_2, \gamma_3),
\label{ac}
\end{equation}
where $f_{ijk}$ are some noncommutative polynomials in four 
variables.\

For every action $\cdot$ there exist an invertible matrix
$M=\left(
\begin{array}{cc}
m_{11} & m_{12} \\
m_{21} & m_{22}
\end{array}
\right)$
$\in {\it C}(1,3)_{2\times 2}$, such that
$$
c_{ij}\cdot v=\sum m_{ik}vm^*_{kj},
$$
where 
$
\left(
\begin{array}{cc}
m^*_{11} & m^*_{12} \\
m^*_{11} & m^*_{12}
\end{array}
\right)
$=$M^{-1}$ (see Skolem-Noether theorem for Hopf algebras
\cite{Ko}\cite{Mo}). The action $\cdot$ is called inner if the map
$c_{ij}\rightarrow m_{ij}$ defines an algebra homorphism
$\varphi: GL_2\rightarrow {\it C}(1,3)$. Since the algebra
${\it C}(1,3)$ is isomorphic to the algebra of $4\times 4$
matrices, the homorphism ${\it C}(1,3)$ defines (and is
defined by) a four dimensional module over (the algebraic
structure of) $GL_2$, or, equivalently, a four dimensional
representation of $GL_2$.\

If $\varphi(c_{12}c_{21})$=$0$, then by definition in Figure 1 
the representation $\varphi$ is given for an essentially more 
simple structure, generated by two commuting ``quantum spinors" 
$(c_{21},c_{11})$ and $(c_{22},c_{12})$. In the case when 
$\varphi(c_{12}c_{21})$$\neq 0$ we say that the inner action 
defined by $\varphi$ has nonzero {\it perturbation}.\

If we add the formal inverse $c^{-1}_{11}$, then the algebraic structure 
of Dipper-Donkin quantization
$GL_2$ is generated by the elements in the following diagram.\

\hspace*{2cm}
\begin{picture}(0,80)
\put(115.29,2.39){\makebox(0,0)[cc]{$c_{21}$}}
\put(168.32,2.39){$a_{22}=c_{11}^{-1}d$}
\put(168.32,50.46){$a_{12}=c_{11}^{-1}c_{12}$}
\put(115.29,50.46){\makebox(0,0)[cc]{$c_{11}$}}
\put(117.50,7){\line(1,1){40}}
\put(122,46){\line(1,-1){40}}
\put(115.5,7){\vector(0,1){40}}
\put(168.5,7){\vector(0,1){40}}
\put(120.5,3){\line(1,0){40}}
\put(120.5,51.0){\line(1,0){40}}
\end{picture}
\vspace*{.5cm}
\label{fig1}
\newline 
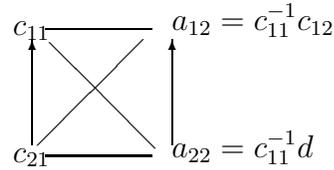
\centerline{Figure 2. Diagramatic representation of}
\centerline{ Dipper-Donkin algebra. 2.}\\

\vspace{1cm}

From here, it follows straightforward that, up to invertibility 
of $c_{11}$, the algebraic structure of
$GL_2$ can be considered like a tensor product 
${\aleph}\otimes {\aleph}$ where ${\aleph}$ is
the quantum plane.\

We say that the representation of the $q$-spinor $xy$=$qyx$,
$x\rightarrow A$, $y\rightarrow B$ is {\it admissible} if 
there exists $C$ such that $x\rightarrow C$, 
$y\rightarrow B$ and $x\rightarrow C$, $y\rightarrow A$ are also 
representations of $q$-spinor with $CB\neq 0$. In other words it
means that $d\rightarrow A$, $c_{12}\rightarrow B$, 
$c_{21}\rightarrow C$ is a representation of the subalgebra of $GL_2$, 
generated by $d$, $c_{12}$, $c_{21}$ with $CB\neq 0$.\

Following the method already developed for studying the actions of 
$GL_q(2,C)$ on the Clifford algebra ${\it C}(1,3)$ \cite{aqg} we can 
construct all inner action of Dipper-Donkin quantization on 
this Clifford algebra, provided $q^m\neq 1$. 
{\it We also can provide with the corresponding operator algebra 
$\Re$ (namely the image of the representation), the algebra of 
invariants $I$ which is equal to the centralizer of $\Re$ in 
${\it C}(1,3)$, and the {\it perturbation} of the representation
\cite{suemi}}. We define $c_{ij}\rightarrow C_{ij}$ to be a finite 
dimensional representation of the quantum $GL_2$.\

Here, we summarize the method used in \cite{suemi} and \cite{aqg} to 
find and classify all possible inner actions of $GL_q(2,C)$ and
$GL_2$ on ${\it C}(1,3)$.  The Hopf algebra $GL_q(2,C)$ is made out of
$q$-spinors, classical spinors and a kind of perturbed spinors. The first
and second types of spinors are defined as it is shown in this paper. For 
the perturbed spinors; in $GL_q(2,C)$, we consider 
$p=(q-q^{-1})a_{12}a_{21}$ meanwhile $p=(q-1)c_{12}c_{21}$ for $GL_2$. 
Actually, $GL_q(2,C)$ can be presented as follows\
\vspace*{1.0cm}
\begin{equation} 
\begin{picture}(80,80)
\put(15.29,2.39){\makebox(0,0)[cc]{$a_{21}$}}
\put(60.0,2.59){\makebox(0,0)[cc]{$a_{22}$}}
\put(60.0,45.46){\makebox(0,0)[cc]{$a_{12}$}}
\put(15.29,45.46){\makebox(0,0)[cc]{$a_{11}$}}
\put(22,45.46){\vector(1,0){28}} 
\put(22,2.39){\vector(1,0){28}}
\put(15.29,38){\vector(0,-1){28}}
\put(58.32,38){\vector(0,-1){28}}
\put(20.0,8.0){\line(1,1){30}}
\put(25.0,33.0){$\cdot$}
\put(30.0,28.0){$\cdot$}
\put(35.0,23.0){$\cdot$}
\put(40.0,18.0){$\cdot$}
\put(45.0,13.0){$\cdot$}
\put(50.0,8.0){$\cdot$}
\end{picture}
\label{p1}
\end{equation}
\newline 
\centerline
{Figure 3. Diagramatic representation of $GL_q(2,C)$.}\\

\vspace{1cm}
Being $a_{11}$, $a_{12}$, $a_{21}$, $a_{22}$, $d^{-1}$ generators of the
algebra. Here $d$=$a_{11}a_{22}-qa_{12}a_{21}$. By arrows we again 
denote $q$-spinors, by the straight line commuting elements (classical 
spinors) and by dots, perturbed spinors.\

For any Hopf algebra that can be presented in terms of combinations of 
$q$-spinors, classical spinors and perturbed spinors (whatever the 
definition of perturbed spinors is), we can 
use our method to find and fully classify all inner actions on
any given algebra. So far we have used ${\it C}(1,3)$ as a particular
example on which $GL_q(2,C)$ and $GL_2$ are acting but we want to stress 
that our method allows us to use any other algebra. We choose 
${\it C}(1,3)$ since we want to have a realization isomorphic 
to the algebra of $4\times 4$ complex matrices, having in mind further
applications in quantum field theory. Moreover, in this paper some results 
are given in terms of matrix units.\

At first, we study all possible $q$-spinor representations on 
${\it C}(1,3)$ such that $q^3,q^4\neq 1$ and analize the equivalence
of representations. Then, we simplify the algebraic structure of
$GL_q(2,C)$ by defining an auxiliary algebra. Form this, we find the
representation of $SL_q(2,C)$ (provided $q^m\neq 1$) and {\it connected}
$GL_q(2,C)$. \cite{aqg}\

Our classification scheme uses this {\it connection} and follows
straightforward \cite{aqg}. The operator algebras, the quantum 
determinants and the invariants of the corresponding inner actions can
also be presented. Two of these given representations are equivalent if
and only if they are equal to each other. Form this, we learn that, for
$GL_q(2,C)$, the quantum determinants are the only quantum invariants.\

The Dipper-Donkin algebra $GL_2$ is also generated by $q$-spinors, 
classical spinors and a perturbated spinor, like $GL_q(2,C)$. Our method 
can be applied with some extra conditions
coming from the particular structure of $GL_2$ \cite{suemi}. Following 
the steps reported above we find, for  $GL_2$, that the corresponding 
algebra of invariants equals the centralizer of coefficients of $M$. By 
definition, the centralizer commutes with all elements of the algebra, 
thereby defining Hamiltonians with conserved energy.\ 

Whenever the algebra of invariants
is in the field (complex for our case) the Hamiltonian obtained is
trivial. Besides, if and only if $I\in \Re$ the corresponding Hamiltonian 
can be defined, since this can be writen in terms of the generators for
$GL_2$. In any other case the invariant algebra cannot be used
to construct Hamiltonias of a quantum chain.

\vspace{1cm}

\section{Quantum chains with global Dipper-Donkin symmetry.}        
In this Section we learn about global quantum group symmetry. Although 
this does not lead straightforward to integrability, we hope to be able 
of giving some additional information related to this subject.\

A method to construct quantum chains with symmetry associated 
to the algebra of functions on a particular quantum group 
has been recently used, for $q$ being root of unity \cite{Ali}.\ 

In this section we address the role of the Dipper-Donkin quantum group
as a global symmetry and show that, contrary with general believe
( see e.g. \cite{Ali}), the mixing of the generators in the coproduct
of the corresponding coalgebra of a quantum group is not sufficient 
condition to construct a non trivial Hamiltonian for quantum chains. See,
for instance, the constructions shown in references \cite{Ri}, \cite{R1}, 
\cite{R2} and \cite{R3}.\

A quantum chain with global quantum group symmetry can be defined as 
follows  \cite{Ali}, \cite{Ri} ; to each site $j$$=1,...,L$ of the chain, 
we assign a representation $\pi_j$. We write the Hamiltonian 
$$
H=\sum^{L-1}_{j=1} id\otimes...\otimes id\otimes 
H_j\otimes id\otimes...\otimes id,
$$
where $H_j$ acts on sites $j$ and $j+1$ as
$$
H_j=\left (\pi_j\otimes \pi_{j+1}\right)[Q_j(\Delta(C))].
$$
Here $j$ denotes the site of the lattice, $C$ is a central
element of the algebra, by $\pi_j$ we mean a representation 
for the algebra and $Q_j$ is a polynomial 
function.\

Let us study, at first, when $q$ is root of unity. In this case, 
this polynomial function can be taken of degree $d\leq p$ 
where the integer $p$ is characterized by the value of $q$ ($q^p$=$1$) as
is done in references \cite{Ali} and \cite{Ri}.\

If $q$ is a root of unity ($q^p=1$), all the elements $c^p_{ij}$ are 
central. In this case, one can uniquely define a state $|0\rangle$ which 
is a common eigenvector of $c_{22}$ and $c_{21}$ with eigenvalues $\alpha$ 
and $\alpha\beta$ ($\alpha$ and $\beta$ being arbitrary constants), 
respectively. Then we build the space $V$ as the linear span of the vectors 
$|n\rangle=c^n_{12}|0\rangle$, $0\leq n\leq p-1$. We can show that $V$ is 
an invariant vector space under the action of the Dipper-Donkin quantum 
group. Thus, we construct $\pi_j$ as follows\
\begin{equation}
\vbox{\halign{# \hfil &\ # \hfil \cr
$c_{12}|n\rangle=|n+1\rangle \mbox{   for $n<p-1$\hspace*{1.7cm}}$& 
$c_{22}|n\rangle=\alpha|n\rangle$\cr 
$c_{11}|n\rangle=\beta|n+1\rangle \mbox{   for $n<p-1$\hspace*{1.7cm}}$& 
$c_{12}|p-1\rangle=\eta |0\rangle$\cr
$c_{21}|n\rangle=q^n\alpha \beta|n\rangle$&  
$c_{11}|p-1\rangle=\beta \eta|0\rangle.\nonumber$\cr }}
\end{equation}

Here $\eta$ is the central value of $c_{12}^p$. All the parameters are 
independent.\

We would like to remark that, in spite of the mixing of the generators in 
the coproduct for the Dipper-Donkin algebra which has the same
structure than the coproduct defined for $GL_q(n)$, the Dipper-Donkin 
algebra leads only to trivial Hamiltonians (proportional to the
identity). This is true, even for the two states quantum chains for which 
it is known that a $GL_q(2)$ global symmetry can be implemented 
\cite{Ali}.\

Now, let us study the case $q^m\neq 1$. We introduce here an alternative way 
to build up Hamiltonians with Dipper-Donkin quantum global symmetry. Since 
for any Dipper-Donkin quantum group the quantum determinant is group-like 
but not central and the invariants $I$ are central, we can 
define a Hamiltonian as follows
$$
H_j=(\pi_j\otimes \pi_{j+1})[Q_j(\Delta (I)].
$$
Our method works only for $q^m\neq 1$. Here $Q_j$ is any polinomial 
function.\

We propose, as a particular case, to study the Dipper-Donkin 
algebra like a global symmetry of four states quantum chains. This is done 
by searching all possible finite dimensional representations of the 
Dipper-Donkin group on the algebra of $4\times 4$ complex matrices on 
which a well defined coproduct for the algebra of invariants can be 
applied. Here $\pi_j$ is one of a this representations.\

In next Table we give the full set of all possible inner actions
that, being non trivial, are in the operator algebra $\Re$; thereby 
properly defining $\Delta (I)$ and a corresponding Hamiltonian. Each 
particular case provides a Hamiltonian with quantum Dipper-Donkin global 
symmetry for a four states quantum chain. An important result of this
paper is that for all the cases reported in next Table
{\it the Hamiltonians for four states quantum chains with Dipper-Donkin
global symmetry have the unique form 
$H_j=(\pi_j\otimes \pi_{j+1})[Q_j(\Delta (A_jd+B_jC_{11}+C_jC_{22})]$} being
$A_j$, $B_j$ and $C_j$ some constants also given in the Table. This, 
together with the representation of the $GL_2$ 
generators, straightforward leads to a Hamiltonian writen in matrix units, 
Dirac gamma matrices or ``mass" ($m_{\pm}$=$\left(1\pm \gamma_0\right)/2$)
and ``spin" ($s_{\uparrow/\downarrow}$=$\left(1\pm i\gamma_{12}\right)/2$) 
operators.\

Let us now introduce some concrete examples for 
each case (namely particular form of the quantum determinant). To reach 
this goal we obtain at first $(\pi_j\otimes \pi_{j+1})[Q_j(\Delta (C_{11})]$. 
 Explicitly,
$$
(\pi_j\otimes \pi_{j+1})[Q_j(\Delta (C_{11}\otimes C_{11})+
Q_j(\Delta (C_{12}\otimes C_{21})]
=$$
$$
\pi_j[Q_j(C_{11}]\otimes 
\pi_{j+1}[Q_{j}(C_{11}]+\pi_j[Q_j(C_{12}]\otimes 
\pi_{j+1}[Q_{j}(C_{21}].
$$ 
We consider the simplest case and take $Q_j$ to be lineal.\

In a similar way we get $(\pi_j\otimes \pi_{j+1})[Q_j(\Delta (C_{22})]$. 
Explicitly 
$$
(\pi_j\otimes \pi_{j+1})[Q_j(\Delta (C_{21}\otimes C_{12})+
Q_j(\Delta (C_{22}\otimes C_{22})]
=$$
$$
\pi_j[Q_j(C_{21}]\otimes 
\pi_{j+1}[Q_{j}(C_{12}]+\pi_j[Q_j(C_{22}]\otimes 
\pi_{j+1}[Q_{j}(C_{22 }].
$$ 
Again we consider the simplest case and
take $Q_j$ to be lineal. At last, we obtain 
$(\pi_j\otimes \pi_{j+1})[Q_j(\Delta (d)]$. Explicitly
$$
(\pi_j\otimes \pi_{j+1})[Q_j(\Delta (d\otimes d)]
=
$$
$$
\pi_j[Q_j(C_{11}C_{22})]\otimes \pi_{j+1}[Q_{j}(C_{11}C_{22})]-
\pi_j[Q_j(C_{11}C_{22})]\otimes \pi_{j+1}[Q_{j}(C_{12}C_{21})]-
$$
$$
\pi_j[Q_j(C_{12}C_{21})]\otimes \pi_{j+1}[Q_{j}(C_{11}C_{22})]+
\pi_j[Q_j(C_{12}C_{21})]\otimes \pi_{j+1}[Q_{j}(C_{12}C_{21})].
$$ 
As usual, we consider $Q_j$ to be lineal.\

We are now ready to present concrete examples.\newline
{\bf CASE 2.2)}
$$
H_j=\left(m_+s_{\uparrow}+\alpha_jm_+s_{\downarrow}+\beta_jm_-s_{\uparrow}
+\gamma_jm_-s_{\downarrow}\right)\otimes 
$$
$$
\left(m_+s_{\uparrow}+
\alpha_{j+1}m_+s_{\downarrow}+\beta_{j+1}m_-s_{\uparrow}
+\gamma_{j+1}m_-s_{\downarrow}\right)
$$
{\bf CASE 3.5)}
$$
H_j=\left(q^2/\alpha_jm_+s_{\uparrow}+q^2/\alpha_jm_+s_{\downarrow}+
m_-s_{\uparrow}+m_-s_{\downarrow}-q^2/\alpha_j^2m_+
(\gamma_1+i\gamma_2)\gamma_3/2\right)\otimes 
$$
$$
\left(q^2/\alpha_{j+1}m_+s_{\uparrow}+q^2/\alpha_{j+1}m_+s_{\downarrow}+
m_-s_{\uparrow}+m_-s_{\downarrow}-q^2/\alpha_{j+1}^2m_+
(\gamma_1+i\gamma_2)\gamma_3/2\right)
$$
{\bf CASE 4.4)}
$$
H_j=
$$
$$
A_j\left(\alpha_jm_+s_{\uparrow}+q^2m_+s_{\downarrow}+
qm_-s_{\uparrow}+m_-s_{\downarrow}\right)\otimes
\left(\alpha_{j+1}m_+s_{\uparrow}+q^2m_+s_{\downarrow}+
qm_-s_{\uparrow}+m_-s_{\downarrow}\right)
$$
$$
+B_j\left(\delta_jm_+s_{\uparrow}+qm_+s_{\downarrow}+
qm_-s_{\uparrow}+m_-s_{\downarrow}\right)\otimes
\left(\delta_{j+1}m_+s_{\uparrow}+qm_+s_{\downarrow}+
qm_-s_{\uparrow}+m_-s_{\downarrow}\right)
$$
$$
-B_j\gamma_j\beta_{j+1}m_+(\gamma_1-i\gamma_2)/2\otimes
m_-(-\gamma_1+i\gamma_2)\gamma_3/2+
$$
$$
C_j\left(\alpha_j/\delta_jm_+s_{\uparrow}+qm_+s_{\downarrow}+
m_-s_{\uparrow}+m_-s_{\downarrow}\right)\otimes
\left(\alpha_{j+1}/\delta_{j+1}m_+s_{\uparrow}+qm_+s_{\downarrow}+
qm_-s_{\uparrow}+m_-s_{\downarrow}\right)
$$
$$
-C_j\gamma_{j+1}\beta_{j}m_-(-\gamma_1+i\gamma_2)\gamma_3/2\otimes
m_+(\gamma_1-i\gamma_2)/2.
$$
{\bf CASE 5.5)}
$$
H_j=\left(q^2\beta_jm_+s_{\uparrow}+q^2\beta_jm_+s_{\downarrow}+
q^2\beta_jm_-s_{\uparrow}+\delta_jm_-s_{\downarrow}+\beta_jm_+
(\gamma_1+i\gamma_2)\gamma_3/2\right)\otimes 
$$
$$
\left(q^2\beta_{j+1}m_+s_{\uparrow}+q^2\beta_{j+1}m_+s_{\downarrow}+
q^2\beta_{j+1}m_-s_{\uparrow}+\delta_{j+1}m_-s_{\downarrow}+\beta_{j+1}m_+
(\gamma_1+i\gamma_2)\gamma_3/2\right)
$$
{\bf CASE 6.10)}
$$
H_j=A_j\left(q^2m_+s_{\uparrow}+qm_+s_{\downarrow}+m_-s_{\uparrow}+
m_-s_{\downarrow}+m_-(\gamma_1+i\gamma_2)\gamma_3/2\right)\otimes
$$
$$
\left(q^2m_+s_{\uparrow}+qm_+s_{\downarrow}+m_-s_{\uparrow}+
m_-s_{\downarrow}\right)- 
$$
$$
B_j(1+m_-(\gamma_1+i\gamma_2)\gamma_3/2\otimes
(1+m_-(\gamma_1+i\gamma_2)\gamma_3/2+
$$
$$
C_j\left(q^2m_+s_{\uparrow}+qm_+s_{\downarrow}+m_-s_{\uparrow}+
m_-s_{\downarrow}\right)\otimes
\left(q^2m_+s_{\uparrow}+qm_+s_{\downarrow}+m_-s_{\uparrow}+
m_-s_{\downarrow}\right)
$$
From these concrete examples it is straigthforward to see that the 
Hamiltonians we have introduced are somehow related with Ashkin-Teller
model written in terms of Ising spin; namely with each site $i$ we associate
two spins (in our case, $m_{\pm}$ and $s_{\uparrow/\downarrow}$).\

Some additional standard symmetries can also be identified. For example,
in CASE 3.5), there is an invariance under the transformation
$s_{\uparrow}$$\rightarrow$$s_{\downarrow}$.\

Finally, in next Table, we explicitly show all the inner actions of $GL_2$ 
on ${\it C}(1,3)$ which can provide with non trivial Hamiltonians for four 
states quantum chains with Dipper-Donkin quantum global symmetry, the 
operator algebra $\Re$, the algebra of invariants $I$ and the 
value of the coefficients $A_j$, $B_j$ and $C_j$ in the unique expression 
for the corresponding Hamiltonians. In all the reported cases the 
perturbation is zero.\

\centerline{\Huge  CASE 1) $d$=$diag(q^2,q,1,1)$}
\vspace*{1.5cm}
\vbox{\offinterlineskip
\hrule
\halign{&\vrule#&\strut \quad \hfil# \hfil
 &\vrule#&\quad \hfil# \hfil&\vrule#& \quad# \hfil \cr
height5pt&\omit&& \omit&& \omit&\cr&
        {\bf CASE 1.1)}               
&&
$\begin{array}{c}
C_{12}=\alpha e_{12}+\beta e_{24}\\
C_{21}=0\\
C_{11}={\bf 1}+e_{34}\\
C_{22}=q^2e_{11}+qe_{22}+e_{33}+e_{44}-e_{34}
\end{array}
$
&& 
        $\Re =\left(\matrix{*&*&0&0\cr 0&*&0&*\cr 0&0&
\epsilon &*\cr 0&0&0&\epsilon \cr}
\right)$&\cr
height5pt&\omit&& \omit&&\omit&\cr
\noalign{\hrule }
height5pt&\omit&& \omit&&\omit&\cr    
&
$\matrix{dim \Re \cr 6 \cr dim I \cr 2}$
 &&
$I =\left(\matrix{\alpha&0&0&0\cr 0&\alpha&0&0\cr 0&0&
\alpha & \gamma\cr 0&0&0&\alpha \cr}
\right)$
&& 
         $A_j=C_j=0$ &\cr
height5pt&\omit&& \omit&& \omit&\cr
\noalign{\hrule }
height2pt&\omit&& \omit&& \omit&\cr
\noalign{\hrule }
        &{\bf CASE 1.2)}  
&&
$\begin{array}{l}
C_{12}=0\\
C_{21}=\alpha e_{21}+\beta e_{32}\\
C_{11}=e_{11}+q^{-1}e_{22}+q^{-2}e_{33}+q^{-2}e_{44}+e_{34}\\ 
C_{22}=q^{2}{\bf 1}-q^4e_{34}\end{array}
$
&& 
        $\Re =\left(\matrix{*&0&0&0\cr *&*&0&0\cr 0&*&
\epsilon &*\cr 0&0&0&\epsilon \cr}
\right)$&\cr
height5pt&\omit&& \omit&&\omit&\cr
\noalign{\hrule }
height5pt&\omit&& \omit&&\omit&\cr    
&$\matrix{dim \Re \cr 6 \cr 
dim I \cr 2}$ && 
$I =\left(\matrix{\alpha&0&0&0\cr 0&\alpha&0&0\cr 0&0&
\alpha & \beta\cr 0&0&0&\alpha \cr}
\right)$&& 
$A_j=B_j=0$ &\cr
height5pt&\omit&& \omit&& \omit&\cr
\noalign{\hrule }
height2pt&\omit&& \omit&& \omit&\cr
\noalign{\hrule }

	&{\bf CASE 1.3)}   
&&
$\begin{array}{l}
C_{12}=\alpha e_{12} \\ 
C_{21}=\beta e_{32}\\ 
C_{11}= e_{11}+e_{22}+q^{-1}e_{33}+q^{-1}e_{44}+e_{34}\\
C_{22}= q^2e_{11}+qe_{22}+qe_{33}+qe_{44}-q^2e_{34}
\end{array}
$
&& 
        $\Re =\left(\matrix{*&*&0&0\cr 0&*&0&0\cr 0&*&
\epsilon &*\cr 0&0&0&\epsilon \cr}
\right)$&\cr
height5pt&\omit&& \omit&&\omit&\cr
\noalign{\hrule }
height5pt&\omit&& \omit&&\omit&\cr    
&$\matrix{dim \Re \cr 6 \cr 
dim I \cr 2 \cr }$ && 
$I =\left(\matrix{\alpha&0&0&0\cr 0&\alpha&0&0\cr 0&0&
\alpha & \beta\cr 0&0&0&\alpha \cr}
\right)$
           &&
$\begin{array}{c} 
A_j=-C_j=-\alpha q^{-1}\\
B_j=\alpha 
\end{array}$
&\cr
height5pt&\omit&& \omit&& \omit&\cr
\noalign{\hrule }
height2pt&\omit&& \omit&& \omit&\cr
\noalign{\hrule }}}


\centerline{\Huge  CASE 2) $d$=$diag(q^2,q,q,1)$}
\vspace*{1.5cm}
\vbox{\offinterlineskip
\hrule
\halign{&\vrule#&\strut \quad \hfil# \hfil
 &\vrule#&\quad \hfil# \hfil&\vrule#& \quad# \hfil \cr
height5pt&\omit&& \omit&& \omit&\cr&

{\bf CASE 2.1)}   
&&
$\begin{array}{c}
	C_{12}=q\lambda\delta e_{13}\\
	C_{21}=\delta e_{43}\\
	C_{11}=e_{11}+\alpha e_{22}+e_{33}+q^{-1}e_{44} \\
        C_{22}=q^{2}e_{11}+q\alpha^{-1}e_{22}+qe_{33}+qe_{44} \\
	\alpha\neq 1 \\
\end{array}$
&& 
        $\Re =\left(\matrix{*&0&*&0\cr 0&*&0&0\cr 0&0&*&0\cr 
0&0&*&*\cr}
\right)$&\cr
height5pt&\omit&& \omit&&\omit&\cr
\noalign{\hrule }
height5pt&\omit&& \omit&&\omit&\cr    
&$\matrix{dim \Re \cr 6 \cr 
dim I \cr 2}$ 
         &&
$I =\left(\matrix{\alpha&0&0&0\cr 0&\beta&0&0\cr 0&0&
\alpha & 0\cr 0&0&0&\alpha \cr}
\right)$
&& 
$\begin{array}{c}
A_j=-C_j=\frac{\alpha (q^{-1}-1)}{q-1}\\
B_j=\alpha
\end{array}$ 
&\cr
height5pt&\omit&& \omit&& \omit&\cr
\noalign{\hrule }
height2pt&\omit&& \omit&& \omit&\cr
\noalign{\hrule }

        & {\bf CASE 2.2)}   
&&
$\begin{array}{l}
C_{12}=0 ; \hspace{.5cm}\alpha\neq \beta\neq \gamma\neq 1\\ 
C_{21}=0\\ 
C_{11}= e_{11}+\alpha e_{22}+\beta e_{33}+\gamma e_{44}\\
C_{22}= q^2e_{11}+q\alpha^{-1}e_{22}+
q\beta^{-1}e_{33}+\gamma^{-1}e_{44}
\end{array}
$
&& 
        $\Re =\left(\matrix{*&0&0&0\cr 0&*&0&0\cr 0&0&
*&0\cr 0&0&0&* \cr}
\right)$&\cr
height5pt&\omit&& \omit&&\omit&\cr
\noalign{\hrule }
height5pt&\omit&& \omit&&\omit&\cr    
&$\matrix{dim \Re \cr 4 \cr 
dim I \cr 4}$ && 
$I =\left(\matrix{\alpha&0&0&0\cr 0&\beta&0&0\cr 0&0&
\gamma &0\cr 0&0&0&\delta \cr}
\right)$
           && 
         $A_j=C_j=0$
&\cr
height5pt&\omit&& \omit&& \omit&\cr
\noalign{\hrule }
height2pt&\omit&& \omit&& \omit&\cr
\noalign{\hrule }
height5pt&\omit&& \omit&& \omit&\cr


	&{\bf CASE 2.3)}   
&&
$\begin{array}{l}
C_{12}=0 ; \hspace*{.5cm}\alpha\neq \beta\\ 
C_{21}=0\\ 
C_{11}= e_{11}+\alpha e_{22}+\beta e_{33}+\gamma e_{44}\\
C_{22}= q^2e_{11}+\frac{q}{\alpha}e_{22}+\frac{q}{\beta}e_{33}+
\gamma^{-1}e_{44}
\end{array}
$
&&
$\Re =\left(\matrix{*&0&0& 0\cr 0&*&0&0\cr
0&0&*&0\cr0&0&0&*\cr}\right)$
        &\cr
height5pt&\omit&& \omit&&\omit&\cr
\noalign{\hrule }
height5pt&\omit&& \omit&&\omit&\cr    
&$\matrix{dim \Re \cr 4 \cr 
dim I \cr 4 \cr }$ && 
$I =\left(\matrix{\alpha&0&0&0\cr 0&\beta&0&0\cr 0&0&
\gamma & 0\cr 0&0&0&\delta \cr}\right)$
           && 
         $A_j=C_j=0$		
&\cr
height5pt&\omit&& \omit&& \omit&\cr
\noalign{\hrule }
height2pt&\omit&& \omit&& \omit&\cr
\noalign{\hrule }}}
\vbox{\offinterlineskip
\hrule
\halign{&\vrule#&\strut \quad \hfil# \hfil
 &\vrule#&\quad \hfil# \hfil&\vrule#& \quad# \hfil \cr
height5pt&\omit&& \omit&& \omit&\cr&
{\bf CASE 2.4)}   
&&
$\begin{array}{l}
C_{12}=0\\ 
C_{21}=0\\ 
C_{11}= e_{11}+\alpha e_{22}+\alpha e_{33}+\beta e_{44}+e_{23}\\
C_{22}= q^2e_{11}+\frac{q}{\alpha}e_{22}+\frac{q}{\alpha}e_{33}+
\frac{1}{\beta}e_{44}-\frac{q}{\alpha^2}e_{23}
\end{array}
$
&&
$\Re =\left(\matrix{*&0&0& 0\cr 0&\epsilon&*&0\cr
0&0&\epsilon&0\cr0&0&0&*\cr}\right)$
        &\cr
height5pt&\omit&& \omit&&\omit&\cr
\noalign{\hrule }
height5pt&\omit&& \omit&&\omit&\cr    
&$\matrix{dim \Re \cr 4\cr
dim I \cr 4 \cr }$ && 
$I =\left(\matrix{\alpha&0&0&0\cr 0&\gamma&\delta&0\cr 0&0&
\gamma & 0\cr 0&0&0&\beta \cr}\right)$
           && 
$A_j=B_j=0$ &\cr
height5pt&\omit&& \omit&& \omit&\cr
height5pt&\omit&& \omit&& \omit&\cr
\noalign{\hrule }
height2pt&\omit&& \omit&& \omit&\cr
\noalign{\hrule }
height5pt&\omit&& \omit&& \omit&\cr

	&{\bf CASE 2.5)}   
&&
$\begin{array}{l}
C_{12}=q\lambda\delta e_{13}\\ 
C_{21}=\delta e_{43}\\ 
C_{11}= e_{11}+e_{22}+e_{33}+q^{-1}e_{44}+e_{23}\\
C_{22}= q^2e_{11}+qe_{22}+qe_{33}+qe_{44}-qe_{23}
\end{array}
$
&&
$\Re =\left(\matrix{*&0&*&0\cr 0&\epsilon&*&0\cr
0&0&\epsilon&0\cr 0&0&*&*\cr}\right)$
        &\cr
height5pt&\omit&& \omit&&\omit&\cr
\noalign{\hrule }
height5pt&\omit&& \omit&&\omit&\cr    
&$\matrix{dim \Re \cr 6\cr
dim I \cr 2 \cr }$ && 
$I =\left(\matrix{\alpha&0&0&0\cr 0&\alpha&\beta&0\cr 0&0&
\alpha& 0\cr 0&0&0&\alpha \cr}\right)$
           && 
$\begin{array}{c}
A_j=-C_j=\frac{\alpha (q^{-1}-1)}{q-1}\\
B_j=\alpha
\end{array}$ 
&\cr
height5pt&\omit&& \omit&& \omit&\cr
\noalign{\hrule }
height2pt&\omit&& \omit&& \omit&\cr
\noalign{\hrule }
height5pt&\omit&& \omit&& \omit&\cr

	&{\bf CASE 2.6)}   
&&
$\begin{array}{l}
C_{12}=0\\ 
C_{21}=0\\ 
C_{11}= e_{11}+\alpha e_{22}+\alpha e_{33}+\beta e_{44}+e_{23}\\
C_{22}= q^2e_{11}+\frac{q}{\alpha}e_{22}+\frac{q}{\alpha}e_{33}+
\frac{1}{\beta}e_{44}-\frac{q}{\alpha^2}e_{23}
\end{array}
$
&&
$\Re =\left(\matrix{*&0&0&0\cr 0&\epsilon&*&0\cr
0&0&\epsilon&0\cr0&0&0&*\cr}\right)$
        &\cr
height5pt&\omit&& \omit&&\omit&\cr
\noalign{\hrule }
height5pt&\omit&& \omit&&\omit&\cr    
&$\matrix{dim \Re \cr 4\cr
dim I \cr 3 \cr }$ && 
$I =\left(\matrix{\alpha&0&0&0\cr 0&\beta&\gamma&0\cr 0&0&
\beta& 0\cr 0&0&0&\alpha \cr}\right)$
           && 
$A_j=B_j=0$ &\cr
height5pt&\omit&& \omit&& \omit&\cr
\noalign{\hrule }
height2pt&\omit&& \omit&& \omit&\cr
\noalign{\hrule }}}


\centerline{\Huge  CASE 3) $d$=$diag(q^2,q^2,q,1)$}
\vspace*{1.5cm}
\vbox{\offinterlineskip
\hrule
\halign{&\vrule#&\strut \quad \hfil# \hfil
 &\vrule#&\quad \hfil# \hfil&\vrule#& \quad# \hfil \cr
height5pt&\omit&& \omit&& \omit&\cr&
        {\bf CASE 3.1)}               
&&
$\begin{array}{c}
C_{12}=\alpha e_{13}+\gamma e_{34}\\
C_{21}=0\\
C_{11}={\bf 1}+e_{12}\\
C_{22}=q^2e_{11}+q^2e_{22}+e_{33}+e_{44}
\end{array}
$
&& 
        $\Re =\left(\matrix{\epsilon&*&*&*\cr 
0&\epsilon&0&0\cr 0&0&*&*\cr 0&0&0&*\cr}\right)$&\cr
height5pt&\omit&& \omit&&\omit&\cr
\noalign{\hrule }
height5pt&\omit&& \omit&&\omit&\cr    
&
$\matrix{dim \Re \cr 7 \cr dim I \cr 2}$
 &&
$I =\left(\matrix{\alpha&\beta&0&0\cr 
0&\alpha&0&0\cr 0&0&\alpha&0\cr 0&0&0&\alpha \cr}
\right)$
&&          $A_j=C_j=0$&\cr
height5pt&\omit&& \omit&& \omit&\cr
\noalign{\hrule }
height2pt&\omit&& \omit&& \omit&\cr
\noalign{\hrule }
height5pt&\omit&& \omit&& \omit&\cr

	&{\bf CASE 3.2)}   
&&
$\begin{array}{l}
C_{12}=\gamma e_{34}\\ 
C_{21}=\beta e_{32}\\ 
C_{11}= qe_{11}+qe_{22}+e_{33}+e_{44}+e_{12}\\
C_{22}= qe_{11}+qe_{22}+qe_{33}+e_{44}-e_{12}
\end{array}
$
&&
$\Re =\left(\matrix{\epsilon&*&0& 0\cr 0&\epsilon&0&0\cr
0&*&*&*\cr0&0&0&*\cr}\right)$
        &\cr
height5pt&\omit&& \omit&&\omit&\cr
\noalign{\hrule }
height5pt&\omit&& \omit&&\omit&\cr    
&$\matrix{dim \Re \cr 6\cr
dim I \cr 2 \cr }$ && 
$I =\left(\matrix{\alpha&\beta&0&0\cr 0&\alpha&0&0\cr 0&0&
\alpha& 0\cr 0&0&0&\alpha \cr}
\right)$
           && 
$\begin{array}{c}
A_j=-C_j=-q^{-1}\alpha\\
B_j=\alpha
\end{array}$ 
&\cr
height5pt&\omit&& \omit&& \omit&\cr
\noalign{\hrule }
height2pt&\omit&& \omit&& \omit&\cr
\noalign{\hrule }
height5pt&\omit&& \omit&& \omit&\cr

	&{\bf CASE 3.3)}   
&&
$\begin{array}{l}
C_{12}=\gamma e_{34} \\ 
C_{21}=0 \\ 
C_{11}= \alpha e_{11}+\alpha e_{22}+e_{33}+e_{44}+e_{12}\\
C_{22}= \frac{q^2}{\alpha}e_{11}+\frac{q^2}{\alpha}e_{22}+qe_{33}+
e_{44}-\frac{q^2}{\alpha^2}e_{12}
\end{array}
$
&&
$\Re =\left(\matrix{\epsilon&*&0&0\cr 0&\epsilon&0&0\cr
0&0&*&*\cr0&0&0&*\cr}\right)$
        &\cr
height5pt&\omit&& \omit&&\omit&\cr
\noalign{\hrule }
height2pt&\omit&& \omit&&\omit&\cr    
&$\matrix{dim \Re \cr 5\cr
dim I \cr 3 \cr }$ && 
$I =\left(\matrix{\beta&\gamma&0&0\cr 0&\beta&0&0\cr 0&0&
\alpha & 0\cr 0&0&0&\alpha \cr}\right)$
           && 
         $A_j=C_j=0$&\cr
height5pt&\omit&& \omit&& \omit&\cr
\noalign{\hrule }
height2pt&\omit&& \omit&& \omit&\cr
\noalign{\hrule }}}

\vbox{\offinterlineskip
\hrule
\halign{&\vrule#&\strut \quad \hfil# \hfil
 &\vrule#&\quad \hfil# \hfil&\vrule#& \quad# \hfil \cr
height5pt&\omit&& \omit&& \omit&\cr

	&{\bf CASE 3.4)}   
&&
$\begin{array}{l}
C_{12}=\beta e_{23} \\ 
C_{21}=\gamma e_{43}\\ 
C_{11}= \alpha e_{11}+qe_{22}+qe_{33}+e_{44}\\
C_{22}= \frac{q^2}{\alpha}e_{11}+qe_{22}+e_{33}+e_{44}\\
				\\
\alpha\neq q  
\end{array}
$
&&
$\Re =\left(\matrix{*&0&0&0\cr 
0&*&*&0\cr 0&0&*&0\cr 0&0&*&*\cr}
\right)$
        &\cr
height5pt&\omit&& \omit&&\omit&\cr
\noalign{\hrule }
height5pt&\omit&& \omit&&\omit&\cr    
&$\matrix{dim \Re \cr 6\cr 
dim I \cr 2 \cr }$ &&
$I =\left(\matrix{\alpha&0&0&0\cr 0&\beta&0&0\cr 0&0&
\beta& 0\cr 0&0&0&\beta \cr}\right)$
           && 
$\begin{array}{c}
A_j=-B_j=-q^{-1}\alpha\\
C_j=\alpha
\end{array}$ 
&\cr
height5pt&\omit&& \omit&& \omit&\cr
\noalign{\hrule }
height2pt&\omit&& \omit&& \omit&\cr
\noalign{\hrule }
height5pt&\omit&& \omit&& \omit&\cr

	&{\bf CASE 3.5)}   
&&
$\begin{array}{l}
C_{12}=0\\ 
C_{21}=\gamma e_{43}\\ 
C_{11}= \alpha e_{11}+\alpha e_{22}+qe_{33}+e_{44}+e_{12}\\
C_{22}= \frac{q^2}{\alpha}e_{11}+ \frac{q^2}{\alpha}e_{22}+
e_{33}+e_{44}-\frac{q^2}{\alpha^2}e_{12}
\end{array}
$
&&
$\Re =\left(\matrix{\epsilon&*&0& 0\cr 0&\epsilon&0&0\cr
0&0&*&0\cr0&0&*&*\cr}\right)$
        &\cr
height5pt&\omit&& \omit&&\omit&\cr
\noalign{\hrule }
height5pt&\omit&& \omit&&\omit&\cr    
&$\matrix{dim \Re \cr 5\cr
dim I \cr 3 \cr }$ && 
$I =\left(\matrix{\beta&\varphi&0&0\cr 0&\beta&0&0\cr 0&0&
\alpha& 0\cr 0&0&0&\alpha \cr}
\right)$
           && 
$A_j=B_j=0$ &\cr
height5pt&\omit&& \omit&& \omit&\cr
\noalign{\hrule }
height2pt&\omit&& \omit&& \omit&\cr
\noalign{\hrule }
height5pt&\omit&& \omit&& \omit&\cr
	&{\bf CASE 3.6)}   
&&
$\begin{array}{l}
C_{12}=\alpha e_{13}\\ 
C_{21}=\gamma e_{43}\\ 
C_{11}= qe_{11}+\alpha e_{22}+qe_{33}+e_{44}\\
C_{22}= qe_{11}+\frac{q^2}{\alpha}e_{22}+e_{33}+e_{44}\\
			\\
\alpha\neq q
\end{array}
$
&&
$\Re =\left(\matrix{*&0&*&0\cr 0&*&0&0\cr
0&0&*&0\cr 0&0&*&*\cr}\right)$
        &\cr
height5pt&\omit&& \omit&&\omit&\cr
\noalign{\hrule }
height5pt&\omit&& \omit&&\omit&\cr    
&$\matrix{dim \Re \cr 6\cr
dim I \cr 2 \cr }$ && 
$I =\left(\matrix{\alpha&0&0&0\cr 0&\beta&0&0\cr 0&0&
\alpha & 0\cr 0&0&0&\alpha \cr}\right)$
           &&
$\begin{array}{c}
A_j=-B_j=-q^{-1}\alpha\\
C_j=\alpha
\end{array}$
&\cr
height5pt&\omit&& \omit&& \omit&\cr
\noalign{\hrule }
height2pt&\omit&& \omit&& \omit&\cr
\noalign{\hrule }}}

\vbox{\offinterlineskip
\hrule
\halign{&\vrule#&\strut \quad \hfil# \hfil
 &\vrule#&\quad \hfil# \hfil&\vrule#& \quad# \hfil \cr


	&{\bf CASE 3.7)}   
&&
$\begin{array}{l}
C_{12}=\alpha e_{13} \\ 
C_{21}=\gamma e_{43} \\ 
C_{11}= qe_{11}+qe_{22}+qe_{33}+e_{44}+e_{12}\\
C_{22}= qe_{11}+qe_{22}+e_{33}+e_{44}-e_{12}
\end{array}
$
&&
$\Re =\left(\matrix{\epsilon&*&*&0\cr 0&\epsilon&0&0\cr
0&0&*&*\cr0&0&0&*\cr}\right)$
        &\cr
height5pt&\omit&& \omit&&\omit&\cr
\noalign{\hrule }
height2pt&\omit&& \omit&&\omit&\cr    
&$\matrix{dim \Re \cr 6\cr
dim I \cr 2 \cr }$ && 
$I =\left(\matrix{\alpha&\beta&0&0\cr 0&\alpha&0&0\cr 0&0&
\alpha & 0\cr 0&0&0&\alpha \cr}\right)$
           && 
$\begin{array}{c}
A_j=-C_j=-q^{-1}\alpha\\
B_j=\alpha
\end{array}$ 
 &\cr
height5pt&\omit&& \omit&& \omit&\cr
\noalign{\hrule }
height2pt&\omit&& \omit&& \omit&\cr
\noalign{\hrule }}}
\newpage

\centerline{\Huge  CASE 4) $d$=$diag(\alpha,q^2,q,1)$}
\centerline{}
\centerline{\Large $\alpha\neq 0, q^{-1},1,q,q^2,q^3$}
\vspace*{1.5cm}
\vbox{\offinterlineskip
\hrule
\halign{&\vrule#&\strut \quad \hfil# \hfil
 &\vrule#&\quad \hfil# \hfil&\vrule#& \quad# \hfil \cr
height5pt&\omit&& \omit&& \omit&\cr&
        {\bf CASE 4.1)}   
&&
$\begin{array}{l}
C_{12}=0\\ 
C_{21}=\gamma e_{32}+\beta e_{43} \\
C_{11}=\delta e_{11}+q^2 e_{22}+qe_{33}+e_{44}\\
C_{22}=\frac{\alpha}{\delta}e_{11}+e_{22}+e_{33}+e_{44}\\
\end{array}$
&& 
        $\Re =\left(\matrix{*&0&0&0\cr 0&*&0&0\cr 0&*&*&0
\cr 0&0&*&*\cr}
\right)$&\cr
height5pt&\omit&& \omit&&\omit&\cr
\noalign{\hrule }
height5pt&\omit&& \omit&&\omit&\cr    
&$\matrix{dim \Re \cr 6 \cr 
dim I \cr 2}$ &&
$I =\left(\matrix{\beta&0&0&0\cr 0&\gamma&0&0\cr 0&0&
\gamma & 0\cr 0&0&0&\gamma \cr}
\right)$
       && 
        $A_j=B_j=0$ &\cr
height5pt&\omit&& \omit&& \omit&\cr
\noalign{\hrule }
height2pt&\omit&& \omit&& \omit&\cr
\noalign{\hrule }
height5pt&\omit&& \omit&& \omit&\cr&
        {\bf CASE 4.2)}               
&&
$\begin{array}{c}
C_{12}=\beta e_{23}+\gamma e_{34}\\
C_{21}=0\\
C_{11}=\delta e_{11}+e_{22}+e_{33}+e_{44}\\
C_{22}=\frac{\alpha}{\delta}e_{11}+q^2e_{22}+qe_{33}+e_{44}
\end{array}
$
&& 
        $\Re =\left(\matrix{*&0&0&0\cr 0&*&*&0\cr 0&0&
*&*\cr 0&0&0&* \cr}
\right)$&\cr
height5pt&\omit&& \omit&&\omit&\cr
\noalign{\hrule }
height5pt&\omit&& \omit&&\omit&\cr    
&
$\matrix{dim \Re \cr 6 \cr dim I \cr 2}$
 &&
$I =\left(\matrix{\gamma&0&0&0\cr 0&\beta&0&0\cr 0&0&
\beta & 0\cr 0&0&0&\beta \cr}
\right)$
&& 
                  $A_j=C_j=0$&\cr
height5pt&\omit&& \omit&& \omit&\cr
\noalign{\hrule }
height2pt&\omit&& \omit&& \omit&\cr
\noalign{\hrule }
height5pt&\omit&& \omit&& \omit&\cr&
{\bf CASE 4.3)}   
&& 
$\begin{array}{c}
	C_{12}=\beta e_{34}\\
	C_{21}=\gamma e_{32}\\
	C_{11}=\delta e_{11}+qe_{22}+e_{33}+e_{44}\\
        C_{22}=\frac{\alpha}{\delta}e_{11}+qe_{22}+qe_{33}+e_{44}
\end{array}$
&& 
        $\Re =\left(\matrix{*&0&0&0\cr 0&*&0&0\cr 0&*&
*&*\cr 0&0&0&*\cr}
\right)$&\cr
height5pt&\omit&& \omit&&\omit&\cr
\noalign{\hrule }
height5pt&\omit&& \omit&&\omit&\cr    
&$\matrix{dim \Re \cr 6 \cr 
dim I \cr 2}$ 
         &&
$I =\left(\matrix{\gamma&0&0&0\cr 0&\beta&0&0\cr 0&0&
\beta & 0\cr 0&0&0&\beta \cr}
\right)$
&& 
$\begin{array}{c}
A_j=-C_j=-q^{-1}\beta\\
B_j=\beta
\end{array}$          
&\cr
height5pt&\omit&& \omit&& \omit&\cr
\noalign{\hrule }
height2pt&\omit&& \omit&& \omit&\cr
\noalign{\hrule }}}
\vbox{\offinterlineskip
\hrule
\halign{&\vrule#&\strut \quad \hfil# \hfil
 &\vrule#&\quad \hfil# \hfil&\vrule#& \quad# \hfil \cr
height5pt&\omit&& \omit&& \omit&\cr
        &{\bf CASE 4.4)}  
&&
$\begin{array}{l}
C_{12}=\gamma e_{23}\\
C_{21}=\beta e_{43}\\
C_{11}=\delta e_{11}+qe_{22}+qe_{33}+e_{44}\\ 
C_{22}=\frac{\alpha}{\delta}e_{11}+qe_{22}+e_{33}+e_{44}
\end{array}
$
&& 
        $\Re =\left(\matrix{*&0&0&0\cr 0&*&*&0\cr 0&0&
*&0\cr 0&0&*&*\cr}
\right)$&\cr
height5pt&\omit&& \omit&&\omit&\cr
\noalign{\hrule }
height5pt&\omit&& \omit&&\omit&\cr    
&$\matrix{dim \Re \cr 6 \cr 
dim I \cr 2}$ && 
$I =\left(\matrix{\gamma&0&0&0\cr 0&\beta&0&0\cr 0&0&
\beta & 0\cr 0&0&0&\beta \cr}
\right)$
&& 
$\begin{array}{c}
A_j=-B_j=-q^{-1}\beta\\
C_j=\beta
\end{array}$ 
 &\cr
height5pt&\omit&& \omit&& \omit&\cr
\noalign{\hrule }
height2pt&\omit&& \omit&& \omit&\cr
\noalign{\hrule }}}
\centerline{\Huge  CASE 5) $d$=$diag(q^2,q^2,q,1)+e_{12}$}
\vspace*{1.5cm}
\vbox{\offinterlineskip
\hrule
\halign{&\vrule#&\strut \quad \hfil# \hfil
 &\vrule#&\quad \hfil# \hfil&\vrule#& \quad# \hfil \cr
height5pt&\omit&& \omit&& \omit&\cr&
        {\bf CASE 5.1)}   
&&
$\begin{array}{l}
C_{12}=0\\ 
C_{21}=\alpha e_{43}\\
C_{11}=q^2\beta e_{11}+q^2\beta e_{22}+q\gamma e_{33}+\gamma e_{44}+
\beta e_{12} \\
C_{22}=\beta^{-1}e_{11}+\beta^{-1}e_{22}+q\gamma^{-1}e_{33}+
\gamma^{-1}e_{44}\\
\end{array}$
&& 
        $\Re =\left(\matrix{\epsilon&*&0&0\cr 0&\epsilon&0&0\cr 0&0&
*&0\cr 0&0&*&*\cr}
\right)$&\cr
height5pt&\omit&& \omit&&\omit&\cr
\noalign{\hrule }
height5pt&\omit&& \omit&&\omit&\cr    
&$\matrix{dim \Re \cr 5 \cr 
dim I \cr 3}$ &&
$I =\left(\matrix{\beta&\gamma&0&0\cr 0&\beta&0&0\cr 0&0&
\alpha & 0\cr 0&0&0&\alpha \cr}
\right)$         
&& 
$A_j=-\gamma B_j-\gamma^{-1}C_j$          
 &\cr
height5pt&\omit&& \omit&& \omit&\cr
\noalign{\hrule }
height2pt&\omit&& \omit&& \omit&\cr
\noalign{\hrule }
height5pt&\omit&& \omit&& \omit&\cr&
        {\bf CASE 5.2)}               
&&
$\begin{array}{c}
C_{12}=\alpha e_{13}+\beta e_{34}\\
C_{21}=0\\
C_{11}=q^2e_{11}+q^2e_{22}+q^2e_{33}+q^2e_{44}+e_{12}\\
C_{22}=e_{11}+e_{22}+q^{-1}e_{33}+q^{-2}e_{44}
\end{array}
$
&& 
        $\Re =\left(\matrix{\epsilon&*&*&0\cr 0&\epsilon&0&0\cr 0&0&
* &*\cr 0&0&0&* \cr}
\right)$&\cr
height5pt&\omit&& \omit&&\omit&\cr
\noalign{\hrule }
height5pt&\omit&& \omit&&\omit&\cr    
&
$\matrix{dim \Re \cr 6 \cr dim I \cr 2}$
 &&
$I =\left(\matrix{\alpha&\beta&0&0\cr 0&\alpha&0&0\cr 0&0&
\alpha & 0\cr 0&0&0&\alpha \cr}
\right)$
&& 
                 $A_j=C_j=0$&\cr
height5pt&\omit&& \omit&& \omit&\cr
\noalign{\hrule }
height2pt&\omit&& \omit&& \omit&\cr
\noalign{\hrule }
height5pt&\omit&& \omit&& \omit&\cr&
{\bf CASE 5.3)}   
&&
$\begin{array}{c}
	C_{12}=\alpha e_{34}\\
	C_{21}=\beta e_{32}\\
	C_{11}=q^2e_{11}+q^{2}e_{22}+qe_{33}+qe_{44}+e_{12}\\
        C_{22}=e_{11}+e_{22}+e_{33}+q^{-1}e_{44}
\end{array}$
&& 
        $\Re =\left(\matrix{\epsilon&*&0&0\cr 0&\epsilon&0&0\cr 0&*&
*&0\cr 0&0&*&* \cr}
\right)$&\cr
height5pt&\omit&& \omit&&\omit&\cr
\noalign{\hrule }
height5pt&\omit&& \omit&&\omit&\cr    
&$\matrix{dim \Re \cr 6 \cr 
dim I \cr 2}$ 
         &&
$I =\left(\matrix{\alpha&\beta&0&0\cr 0&\alpha&0&0\cr 0&0&
\alpha & 0\cr 0&0&0&\alpha \cr}
\right)$
&& 
$\begin{array}{c}
A_j=-B_j=-q^{-1}\beta\\
C_j=\beta
\end{array}$         
&\cr
height5pt&\omit&& \omit&& \omit&\cr
\noalign{\hrule }
height2pt&\omit&& \omit&& \omit&\cr
\noalign{\hrule }}}
\vbox{\offinterlineskip
\hrule
\halign{&\vrule#&\strut \quad \hfil# \hfil
 &\vrule#&\quad \hfil# \hfil&\vrule#& \quad# \hfil \cr
height5pt&\omit&& \omit&& \omit&\cr
        &{\bf CASE 5.4)}  
&&
$\begin{array}{l}
C_{12}=\alpha e_{34}\\
C_{21}=0\\
C_{11}=q^2\beta e_{11}+q^2\beta e_{22}+\gamma e_{33}+\gamma e_{44}
+\beta e_{12}\\ 
C_{22}=\beta^{-1}e_{11}+\beta^{-1}e_{22}+q\gamma^{-1}e_{33}+
\gamma^{-1}e_{44}
\end{array}
$
&& 
        $\Re =\left(\matrix{\epsilon&*&0&0\cr 0&\epsilon&0&0\cr 0&0&
*&*\cr 0&0&0&*\cr}
\right)$&\cr
height5pt&\omit&& \omit&&\omit&\cr
\noalign{\hrule }
height5pt&\omit&& \omit&&\omit&\cr    
&$\matrix{dim \Re \cr 5 \cr 
dim I \cr 3}$ && 
$I =\left(\matrix{\alpha&\gamma&0&0\cr 0&\alpha&0&0\cr 0&0&
\beta & 0\cr 0&0&0&\beta \cr}
\right)$&& 
         $A_j=C_j=0$&\cr
height5pt&\omit&& \omit&& \omit&\cr
\noalign{\hrule }
height2pt&\omit&& \omit&& \omit&\cr
\noalign{\hrule }
height5pt&\omit&& \omit&& \omit&\cr
        & {\bf CASE 5.5)}   
&&
$\begin{array}{c}
	C_{12}=\alpha e_{13}\\
	C_{21}=0\\
	C_{11}=q^2\beta e_{11}+q^2\beta e_{22}+q^2\beta e_{33}
+\delta e_{44}+\beta e_{12}\\
        C_{22}=\beta^{-1}e_{11}+\beta^{-1}e_{22}+q^{-1}\beta^{-1}e_{33}+
\delta^{-1}e_{44}
\end{array}$
&& 
        $\Re =\left(\matrix{\epsilon&*&*&0\cr 0
&\epsilon&0&0\cr 0&0&
*&0\cr 0&0&0&*\cr}
\right)$&\cr
height5pt&\omit&& \omit&&\omit&\cr
\noalign{\hrule }
height5pt&\omit&& \omit&&\omit&\cr    
&$\matrix{dim \Re \cr 5 \cr 
dim I \cr 3}$ && 
$I =\left(\matrix{\alpha&\gamma&0&0\cr 0&\alpha&0&0\cr 0&0&
\alpha& 0\cr 0&0&0&\beta \cr}
\right)$
           && 
         $A_j=C_j=0$ &\cr
height5pt&\omit&& \omit&& \omit&\cr
\noalign{\hrule }
height2pt&\omit&& \omit&& \omit&\cr
\noalign{\hrule }
height5pt&\omit&& \omit&& \omit&\cr
        & {\bf CASE 5.6)}   
&&
$\begin{array}{l}
C_{12}=0 \\ 
C_{21}=\alpha e_{32} \\ 
C_{11}= q^2\beta e_{11}+q^2\beta e_{22}+q\beta e_{33}+
\delta e_{44}+\beta e_{12}\\
C_{22}= \beta^{-1}e_{11}+\beta^{-1}e_{22}+\beta^{-1}e_{33}
+\delta^{-1}e_{44}
\end{array}
$
&& 
        $\Re =\left(\matrix{\epsilon&*&0&0\cr 0&\epsilon&0&0\cr 0&*&
*&0\cr 0&0&0&*\cr}
\right)$&\cr
height5pt&\omit&& \omit&&\omit&\cr
\noalign{\hrule }
height5pt&\omit&& \omit&&\omit&\cr    
&$\matrix{dim \Re \cr 5 \cr 
dim I \cr 3}$ && 
$I =\left(\matrix{\alpha&\gamma&0&0\cr 0&\alpha&0&0\cr 0&0&
\alpha & 0\cr 0&0&0&\beta \cr}
\right)$
           && 
$A_j=-\beta B_j=-b\beta$
 &\cr
height5pt&\omit&& \omit&& \omit&\cr
\noalign{\hrule }
height2pt&\omit&& \omit&& \omit&\cr
\noalign{\hrule }
height5pt&\omit&& \omit&& \omit&\cr
        & {\bf CASE 5.7)}   
&&
$\begin{array}{l}
C_{12}=0 \\ 
C_{21}=0\\ 
C_{11}= q^2\alpha e_{11}+q^2\alpha e_{22}+\beta e_{33}+\gamma e_{44}
+\alpha e_{12}\\
C_{22}= \alpha^{-1}e_{11}+\alpha^{-1}e_{22}+q\beta^{-1}e_{33}+
\gamma^{-1}e_{44}
\end{array}
$
&& 
        $\Re =\left(\matrix{\epsilon&*&0&0\cr 0&\epsilon&0&0\cr 0&0&
*&0\cr 0&0&0&*\cr}
\right)$&\cr
height5pt&\omit&& \omit&&\omit&\cr
\noalign{\hrule }
height5pt&\omit&& \omit&&\omit&\cr    
&$\matrix{dim \Re \cr 4 \cr 
dim I \cr 4}$ && 
$I =\left(\matrix{\gamma&\varphi&0&0\cr 0&\gamma&0&0\cr 0&0&
\beta & 0\cr 0&0&0&\alpha \cr}
\right)$
           && 
         $A_j=C_j=0$&\cr
height5pt&\omit&& \omit&& \omit&\cr
\noalign{\hrule }
height2pt&\omit&& \omit&& \omit&\cr
\noalign{\hrule }}}


\centerline{\Huge  CASE 6) $d$=$diag(q^2,q,1,1)+e_{34}$}
\vspace*{1.5cm}
\vbox{\offinterlineskip
\hrule
\halign{&\vrule#&\strut \quad \hfil# \hfil
 &\vrule#&\quad \hfil# \hfil&\vrule#& \quad# \hfil \cr
height5pt&\omit&& \omit&& \omit&\cr
        & {\bf CASE 6.1)}   
&&
$\begin{array}{c}
	C_{12}=0\\
	C_{21}=\alpha e_{21}\\
	C_{11}=q\beta e_{11}+\beta e_{22}+\gamma e_{33}+\delta e_{44}
+\delta e_{34}\\
        C_{22}=q\beta^{-1}e_{11}+q\beta^{-1}e_{22}+\gamma^{-1}e_{33}+
\delta^{-1}e_{44}\\
\mbox{either }\gamma=\delta=\beta\;\mbox{ or}\;\gamma=\delta=q\beta\\
\mbox{ or }\gamma=\delta\neq \beta\;\mbox{ or }\;\gamma=\delta\neq q\beta  
\end{array}$
&& 
        $\Re =\left(\matrix{*&0&0&0\cr *
&*&0&0\cr 0&0&
\epsilon&*\cr 0&0&0&\epsilon\cr}
\right)$&\cr
height5pt&\omit&& \omit&&\omit&\cr
\noalign{\hrule }
height5pt&\omit&& \omit&&\omit&\cr    
&$\matrix{dim \Re \cr 5 \cr 
dim I \cr 3}$ && 
$I =\left(\matrix{\alpha&0&0&0\cr 0&\alpha&0&0\cr 
0&0&\beta&\gamma\cr 0&0&0&\beta\cr}
\right)$
           && 
$A_j=-q^{-1}\beta B_j$ 
 &\cr
height5pt&\omit&& \omit&& \omit&\cr
\noalign{\hrule }
height2pt&\omit&& \omit&& \omit&\cr
\noalign{\hrule }
height5pt&\omit&& \omit&& \omit&\cr

        & {\bf CASE 6.2)}   
&&
$\begin{array}{l}
C_{12}=0\\ 
C_{21}=\alpha e_{21}\\ 
C_{11}= q\beta e_{11}+\beta e_{22}+\gamma e_{33}+\delta e_{44}+\delta e_{34}\\
C_{22}= q\beta^{-1}e_{11}+q\beta^{-1}e_{22}+\gamma^{-1}e_{33}+
\delta^{-1}e_{44}\\
\mbox{neither }\;\; \gamma=\delta=\beta\;\;\mbox{ nor }\;\;\gamma=\delta=q\beta\\
\hspace*{.5cm}\mbox{nor }\gamma=\delta\neq \beta\;\;\mbox{ nor }\;\;\gamma=\delta\neq q\beta
\end{array}
$
	&& 
        $\Re =\left(\matrix{*&0&0&0\cr *&*&0&0\cr 0&0&
*&*\cr 0&0&0&* \cr}
\right)$&\cr
height5pt&\omit&& \omit&&\omit&\cr
\noalign{\hrule }
height5pt&\omit&& \omit&&\omit&\cr    
&$\matrix{dim \Re \cr 6 \cr 
dim I \cr 2}$ && 
$I =\left(\matrix{\beta&0&0&0\cr 0&\beta&0&0\cr 0&0&
\alpha& 0\cr 0&0&0&\alpha \cr}
\right)$
           && 
$\begin{array}{c}
A_j=\frac{\delta\gamma^{-1}-1}{\gamma-\delta}C_j\\
B_j=\frac{\delta^{-1}-\gamma^{-1}}{\gamma-\delta}C_j
\end{array}$  
 &\cr
height5pt&\omit&& \omit&& \omit&\cr
\noalign{\hrule }
height2pt&\omit&& \omit&& \omit&\cr
\noalign{\hrule }
height5pt&\omit&& \omit&& \omit&\cr
        & {\bf CASE 6.3)}   
&&
$\begin{array}{l}
C_{12}=\alpha e_{24}+\beta e_{12}\\ 
C_{21}=0 \\ 
C_{11}= {\bf 1}+e_{34}\\
C_{22}= q^2e_{11}+qe_{22}+e_{33}+e_{44}
\end{array}
$
&& 
$\Re =\left(\matrix{*&*&0&0\cr 0&*&0&*\cr 0&0&
\epsilon&*\cr 0&0&0&\epsilon\cr}
\right)$
        &\cr
height5pt&\omit&& \omit&&\omit&\cr
\noalign{\hrule }
height5pt&\omit&& \omit&&\omit&\cr    
&$\matrix{dim \Re \cr 6\cr 
dim I \cr 2}$ && 
$I =\left(\matrix{\alpha&0&0&0\cr 0&\alpha&0&0\cr 0&0&
\alpha&\beta\cr 0&0&0&\alpha \cr}
\right)$
           && 
         $A_j=C_j=0$ &\cr
height5pt&\omit&& \omit&& \omit&\cr
\noalign{\hrule }
height2pt&\omit&& \omit&& \omit&\cr
\noalign{\hrule }}}
\vbox{\offinterlineskip
\hrule
\halign{&\vrule#&\strut \quad \hfil# \hfil
 &\vrule#&\quad \hfil# \hfil&\vrule#& \quad# \hfil \cr
height5pt&\omit&& \omit&& \omit&\cr


	&{\bf CASE 6.4)}   
&&
$\begin{array}{l}
C_{12}=\beta e_{12}\\ 
C_{21}=0\\ 
C_{11}= \alpha e_{11}+\alpha e_{22}+\gamma e_{33}+\gamma e_{44}+
\gamma e_{34}\\
C_{22}= \frac{q^2}{\alpha}e_{11}+\frac{q}{\alpha}e_{22}+
\gamma^{-1}e_{33}+\gamma^{-1}e_{44}
\end{array}
$
&&
$\Re =\left(\matrix{*&*&0&0\cr 0&*&0&0\cr 0&0&\epsilon&*\cr
0&0&0&\epsilon}\right)$
        &\cr
height5pt&\omit&& \omit&&\omit&\cr
\noalign{\hrule }
height5pt&\omit&& \omit&&\omit&\cr    
&$\matrix{dim \Re \cr 5\cr
dim I \cr 3 \cr }$ && 
$I =\left(\matrix{\alpha&0&0&0\cr 0&\alpha&0&0\cr 0&0&
\beta&\gamma\cr 0&0&0&\beta \cr}\right)$
          && 
         $A_j=C_j=0$&\cr
height5pt&\omit&& \omit&& \omit&\cr
\noalign{\hrule }
height2pt&\omit&& \omit&& \omit&\cr
\noalign{\hrule }
height5pt&\omit&& \omit&& \omit&\cr


	&{\bf CASE 6.5)}   
&&
$\begin{array}{l}
C_{12}=\beta e_{12}\\ 
C_{21}=0 \\ 
C_{11}= \alpha e_{11}+\alpha e_{22}+\gamma e_{33}+\delta e_{44}+
\delta e_{34}\\
C_{22}= \frac{q^2}{\alpha}e_{11}+\frac{q}{\alpha}e_{22}+\gamma^{-1}e_{33}+
\delta^{-1}e_{44}\\
	\\
\gamma\neq \delta 
\end{array}
$
&&
$\Re =\left(\matrix{*&*&0&0\cr 0&*&0&0\cr
0&0&*&*\cr0&0&0&*\cr}\right)$
        &\cr
height5pt&\omit&& \omit&&\omit&\cr
\noalign{\hrule }
height2pt&\omit&& \omit&&\omit&\cr    
&$\matrix{dim \Re \cr 6\cr
dim I \cr 2 \cr }$ && 
$I =\left(\matrix{\alpha&0&0&0\cr 0&\alpha&0&0\cr 0&0&
\beta& 0\cr 0&0&0&\beta\cr}\right)$
           &&
$\begin{array}{c}
A_j=-\gamma^{-1}\;\;\;C_j=1\\
B_j=\frac{\delta^{-1}-\gamma^{-1}}{\gamma-\delta}
\end{array}$   
 &\cr
height5pt&\omit&& \omit&& \omit&\cr
\noalign{\hrule }
height2pt&\omit&& \omit&& \omit&\cr
\noalign{\hrule }
height5pt&\omit&& \omit&& \omit&\cr

	&{\bf CASE 6.6)}   
&&
$\begin{array}{l}
C_{12}=\delta e_{24}\\ 
C_{21}=0\\ 
C_{11}= \alpha e_{11}+\beta e_{22}+\gamma e_{33}+
\beta e_{44}+\beta e_{34}\\
C_{22}= \frac{q^2}{\alpha}e_{11}+\frac{q}{\beta}e_{22}+
\gamma^{-1}e_{33}+\beta^{-1}e_{44}\\
	\\
\mbox{ either }\alpha=\beta=\gamma\;\mbox { or }\;\gamma^{-1}=q^2\alpha^{-1}=\beta^{-1}
\end{array}
$
&&
$\Re =\left(\matrix{*&0&0&0\cr 0&*&0&*\cr
0&0&\epsilon&*\cr0&0&0&\epsilon\cr}\right)$
        &\cr
height5pt&\omit&& \omit&&\omit&\cr
\noalign{\hrule }
height5pt&\omit&& \omit&&\omit&\cr    
&$\matrix{dim \Re \cr 5\cr
dim I \cr 2 \cr }$ && 
$I =\left(\matrix{\alpha&0&0&0\cr 0&\beta&0&0\cr 0&0&
\beta& 0\cr 0&0&0&\beta \cr}\right)$

           && 
$\begin{array}{c}
B_j=-\beta^{-1}A_j\\
C_j=\frac{1-q}{q\beta^{-1}-\beta^{-1}}A_j\\
\end{array}$   
 &\cr
height5pt&\omit&& \omit&& \omit&\cr
\noalign{\hrule }
height2pt&\omit&& \omit&& \omit&\cr
\noalign{\hrule }}}
\vbox{\offinterlineskip
\hrule
\halign{&\vrule#&\strut \quad \hfil# \hfil
 &\vrule#&\quad \hfil# \hfil&\vrule#& \quad# \hfil \cr
height5pt&\omit&& \omit&& \omit&\cr

	&{\bf CASE 6.7)}   
&&
$\begin{array}{l}
C_{12}=\delta e_{24}\\ 
C_{21}=0\\ 
C_{11}= \alpha e_{11}+\beta e_{22}+\gamma e_{33}+
\beta e_{44}+\beta e_{34}\\
C_{22}= \frac{q^2}{\alpha}e_{11}+\frac{q}{\beta}e_{22}+
\gamma^{-1}e_{33}+\beta^{-1}e_{44}\\
	\\
\mbox{neither }\;\alpha=\beta=\gamma\;\mbox{ nor }\;
1/\gamma=q^2/\alpha=1/\beta
\end{array}
$
&&
$\Re =\left(\matrix{*&0&0&0\cr 0&*&0&*\cr
0&0&*&*\cr0&0&0&*\cr}\right)$
        &\cr
height5pt&\omit&& \omit&&\omit&\cr
\noalign{\hrule }
height5pt&\omit&& \omit&&\omit&\cr    
&$\matrix{dim \Re \cr 6\cr
dim I \cr 2 \cr }$ && 
$I =\left(\matrix{\alpha&0&0&0\cr 0&\beta&0&0\cr 0&0&
\beta& 0\cr 0&0&0&\beta \cr}\right)$
           && 
$\begin{array}{c}
A_j=-\beta B_j\\
B_j=\beta^2B_j
\end{array}$   
 &\cr
height5pt&\omit&& \omit&& \omit&\cr
\noalign{\hrule }
height2pt&\omit&& \omit&& \omit&\cr
\noalign{\hrule }
height5pt&\omit&& \omit&& \omit&\cr

	&{\bf CASE 6.8)}   
&&
$\begin{array}{l}
C_{12}=0\\ 
C_{21}=\beta e_{32}\\ 
C_{11}= \alpha e_{11}+q\gamma e_{22}+\gamma e_{33}+
\delta e_{44}+\delta e_{34}\\
C_{22}= \frac{q^2}{\alpha}e_{11}+\gamma^{-1}e_{22}+
\gamma^{-1}e_{33}+\delta^{-1}e_{44}\\
	\\
\mbox{ either }\alpha=\gamma=\delta\;\mbox{ or }\;\alpha\neq \gamma\neq \delta
\end{array}
$
&&
$\Re =\left(\matrix{*&0&0&0\cr 0&*&0&0\cr
0&*&\epsilon&*\cr0&0&0&\epsilon\cr}\right)$
        &\cr
height5pt&\omit&& \omit&&\omit&\cr
\noalign{\hrule }
height5pt&\omit&& \omit&&\omit&\cr    
&$\matrix{dim \Re \cr 5\cr
dim I \cr 3 \cr }$ && 
$I =\left(\matrix{\alpha&0&0&0\cr 0&\beta&0&0\cr 0&0&
\beta&\gamma\cr 0&0&0&\beta\cr}\right)$
           && 
$\begin{array}{c}
{\mbox For  }\;\;\alpha\neq \gamma\neq \delta\\
A_j=\frac{1-\gamma\delta^{-1}}{\gamma-\delta}\;\;\;C_j=1\\
B_j=\frac{\delta^{-1}-\gamma^{-1}}{\gamma-\delta}\\
{\mbox For  }\;\;\alpha=\gamma=\delta\\
A_j=-\gamma B_j
\end{array}$    &\cr
height5pt&\omit&& \omit&& \omit&\cr
\noalign{\hrule }
height2pt&\omit&& \omit&& \omit&\cr
\noalign{\hrule }
height5pt&\omit&& \omit&& \omit&\cr
	&{\bf CASE 6.9)}   
&&
$\begin{array}{l}
C_{12}=0\\ 
C_{21}=\beta e_{32}\\ 
C_{11}= \alpha e_{11}+q\gamma e_{22}+\gamma e_{33}+
\delta e_{44}+\delta e_{34}\\
C_{22}= \frac{q^2}{\alpha}e_{11}+\gamma^{-1}e_{22}+
\gamma^{-1}e_{33}+\delta^{-1}e_{44}\\
	\\
\mbox{neither }\;\alpha=\gamma=\delta\;\mbox{ nor }\;
\alpha\neq \gamma\neq \delta

\end{array}
$
&&
$\Re =\left(\matrix{*&0&0&0\cr 
0&*&0&0\cr 0&*&*&*\cr 0&0&0&*\cr}
\right)$
        &\cr
height5pt&\omit&& \omit&&\omit&\cr
\noalign{\hrule }
height5pt&\omit&& \omit&&\omit&\cr    
&$\matrix{dim \Re \cr 6\cr 
dim I \cr 2 \cr }$ &&
$I =\left(\matrix{\alpha&0&0&0\cr 0&\beta&0&0\cr 0&0&
\beta& 0\cr 0&0&0&\beta \cr}\right)$
           && 
$\begin{array}{c}
A_j=B_j=0\\
C_j=\frac{\delta^{-1}-\gamma^{-1}}{\gamma-\delta}
\end{array}$  
&\cr
height5pt&\omit&& \omit&& \omit&\cr
\noalign{\hrule }
height2pt&\omit&& \omit&& \omit&\cr
\noalign{\hrule }}}
\vbox{\offinterlineskip
\hrule
\halign{&\vrule#&\strut \quad \hfil# \hfil
 &\vrule#&\quad \hfil# \hfil&\vrule#& \quad# \hfil \cr
height5pt&\omit&& \omit&& \omit&\cr

	&{\bf CASE 6.10)}   
&&
$\begin{array}{l}
C_{12}=0\\ 
C_{21}=0\\ 
C_{11}= {\bf 1}+e_{34}\\
C_{22}= q^2e_{11}+ qe_{22}+e_{33}+e_{44}
\end{array}
$
&&
$\Re =\left(\matrix{*&0&0& 0\cr 0&*&0&0\cr
0&0&\epsilon&*\cr0&0&0&\epsilon\cr}\right)$
        &\cr
height5pt&\omit&& \omit&&\omit&\cr
\noalign{\hrule }
height5pt&\omit&& \omit&&\omit&\cr    
&$\matrix{dim \Re \cr4\cr 
dim I \cr 4 \cr }$ && 
$I =\left(\matrix{\alpha&0&0&0\cr 0&\beta&0&0\cr 0&0&
\gamma& \delta\cr 0&0&0&\gamma \cr}
\right)$
           && 
$A_j$, $B_j$, $C_j$ arbitrary
 &\cr
height5pt&\omit&& \omit&& \omit&\cr
\noalign{\hrule }
height2pt&\omit&& \omit&& \omit&\cr
\noalign{\hrule }
height5pt&\omit&& \omit&& \omit&\cr

	&{\bf CASE 6.11)}   
&&
$\begin{array}{l}
C_{12}=0\\ 
C_{21}=0\\ 
C_{11}= \alpha e_{11}+\beta e_{22}+\gamma e_{33}+
\delta e_{44}+\delta e_{34}\\
C_{22}= \frac{q^2}{\alpha}e_{11}+\frac{q}{\beta}e_{22}+
\gamma^{-1}e_{33}+\epsilon^{-1}e_{44}
    			\\
\alpha\neq \beta\neq \gamma \neq \delta \neq \epsilon
\end{array}
$
&&
$\Re =\left(\matrix{*&0&0&0\cr 0&*&0&0\cr
0&0&*&*\cr0&0&0&*\cr}\right)$
        &\cr
height5pt&\omit&& \omit&&\omit&\cr
\noalign{\hrule }
height5pt&\omit&& \omit&&\omit&\cr    
&$\matrix{dim \Re \cr 5\cr 
dim I \cr 3 \cr }$ && 
$I =\left(\matrix{\alpha&0&0&0\cr 0&\beta&0&0\cr 0&0&
\gamma & 0\cr 0&0&0&\gamma \cr}\right)$
           && 
$\begin{array}{c}
B_j=C_j=\frac{\epsilon^{-1}-\gamma^{-1}}{\gamma-\delta}\\
A_j=\frac{\delta(\gamma^{-1}-\epsilon^{-1})}{\gamma-\delta}
\end{array}$  
 &\cr
height5pt&\omit&& \omit&& \omit&\cr
height5pt&\omit&& \omit&& \omit&\cr
\noalign{\hrule }
height2pt&\omit&& \omit&& \omit&\cr
\noalign{\hrule }}}
\vspace{2cm}
Table. $GL_2$ representations, corresponding operator 
algebras $\Re$, algebra of Invariants which are used to define non trivial 
Hamiltonians for four states quantum chains and
the coefficients in the unique expression for these Hamiltonians are 
presented. The classification used is given by means of the different 
determinants in $GL_2$.

\vspace*{2cm}
\section{Summary and Conclusions.}
We are able of showing the way to construct all possible Hamiltonians
for four states quantum chains with Dipper-Donkin global symmetry, for
$q^m\neq 1$. This
is done, although the Dipper-Donkin quantum group has
non central but group-like determinant. We use the algebra of invariants for
the actions of $GL_2$ on ${\it C}(1,3)$ which corresponds to
the centralizer of the operator algebra, or image of the representation.\

It is straighforward to see that in all the possible cases the 
perturbation of the corresponding action is zero. Moreover, 
there are only few cases where all generators are not null. This 
occurs whenever the algebra of invariants is of the form 
$\alpha {\bf 1}+\beta e_{ii+1}$ or $\alpha e_{ii}+\beta e_{jj}$ (being
$(ii)$=$(11)$ and $(jj)$=$(22)$=$(33)$=$(44)$ or 
$(ii)$=$(22)$ and $(jj)$=$(11)$=$(33)$=$(44)$). We find out that all
the Hamiltonians for four states quantum chains with Dipper-Donkin global
symmetry have the following unique form
$$
H=\sum^{L-1}_{j=1}id\otimes...\otimes id\otimes(\pi_j\otimes \pi_{j+1})
[Q_j(\Delta(A_jd+B_jC_{11}+C_jC_{22}))]\otimes id...\otimes id,
$$
and we provide with the specific values for $A_j$, $B_j$ and $C_j$ in all
cases. Some concrete examples, written in terms of 
$m_{\pm}$ and $s_{\uparrow/\downarrow}$ are also 
introduced.\

We report elsewhere \cite{elec} the complete classification of all 
inner actions of the Dipper-Donkin quantum group on the ${\it C}(1,3)$
algebra. In that paper can be seen how all invariants of the corresponding
CASE 4) generate trivial Hamiltonians (this is not shown in current
paper), all invariants of the corresponding CASE 5) are $C\otimes C$ (this 
is CASE 4) in current paper), and finally all invariants of the corresponding 
CASE 6) (in current paper CASE 5)) are diagonal plus $\beta e_{12}$. It is 
also remarkable that all invariant algebras used to construct Hamiltonians 
with Dipper-Donkin global symmetry for four states quantum chains, are 
either diagonal or have elements in the diagonal plus $\beta e_{ii+1}$.

\section{Acknowledgments}
The author wishes to thank V. Kharchenko for helpful discussion and
CONACYT-M\'exico for partial support under grant No 4336-E.

\end{document}